\documentclass[10pt]{amsart}
\usepackage[matrix,arrow]{xy}
\usepackage{amssymb}
\usepackage{amsmath}
\usepackage{amsfonts}
\usepackage{epsfig}
\usepackage{color}
\usepackage{epstopdf}
\usepackage{graphicx}
\usepackage{amsthm}
\usepackage{enumerate}
\usepackage[mathscr]{eucal}
\usepackage{verbatim}
\usepackage{bm}
\usepackage{mathtools}
\usepackage{graphicx}

\makeatletter
\g@addto@macro\bfseries{\boldmath}
\makeatother

\usepackage[bookmarks=true,hyperindex,pdftex,colorlinks,citecolor=blue,
linkcolor=blue,urlcolor=blue]{hyperref}
\usepackage{tikz}
\usetikzlibrary{matrix,arrows.meta}

\usepackage{orcidlink}

\newtheorem{theorem}{Theorem}[section]
\newtheorem{lemma}[theorem]{Lemma}
\newtheorem{proposition}[theorem]{Proposition}
\newtheorem{corollary}[theorem]{Corollary}

\theoremstyle{remark}
\newtheorem{remark}[theorem]{Remark}
\theoremstyle{definition}

\newcommand{\K}{\mathcal K}

\newcommand{\Fin}{\mathcal{F}}

\newcommand{\rank}{\operatorname{rank}}
\newcommand{\ran}{\operatorname{ran}}

\newcommand{\id}[1]{\mathrm{id}_{#1}}
\newcommand{\ip}[2]{\left\langle #1,#2\right\rangle}
\newcommand{\N}{\mathbb{N}}

\newcommand{\R}{\mathbb{R}}
\newcommand{\Read}{\mathcal{R}}
\newcommand{\sRead}{\widetilde{\mathcal{R}}}
\newcommand{\X}{\mathfrak{X}}
\newcommand{\NA}{\operatorname{NA}}
\newcommand{\Lin}{\mathcal{L}}
\newcommand{\tn}[1]{|\!|\!|#1|\!|\!|}

\newcommand{\supp}{\operatorname{supp}}
\newcommand{\spn}{\operatorname{span}}

\renewcommand{\leq}{\leqslant}
\renewcommand{\geq}{\geqslant}

\title[Norm-attaining operators on renormings of $c_0$]
{Norm-attaining operators on renormings of $c_0$} 

\author[M.~Han]{Manwook Han\orcidlink{0009-0000-7250-1814}}
\address[M.~Han]{Department of Mathematics, Chungbuk National University, Cheongju, Chungbuk 28644, Republic of Korea\newline
	\href{https://orcid.org/0009-0000-7250-1814}{ORCID: \texttt{0009-0000-7250-1814}}}
\email{\texttt{626086@chungbuk.ac.kr}}
\urladdr{\url{https://sites.google.com/view/manwookhan}}

\author[M.~Mart\'in]{Miguel Mart\'in\orcidlink{0000-0003-4502-798X}}
\address[M.~Mart\'in]{Department of Mathematical Analysis and Institute of Mathematics (IMAG), University of Granada, 18071 Granada, Spain \newline \href{http://orcid.org/0000-0003-4502-798X}{ORCID: \texttt{0000-0003-4502-798X} } }
\email{\texttt{mmartins@ugr.es}}
\urladdr{\url{http://www.ugr.es/~mmartins/}}

\author[D. L. Rodríguez-Vidanes]{Daniel L. Rodríguez-Vidanes\orcidlink{0000-0002-1016-096X}}
	\address[D. L. Rodríguez-Vidanes]{Grupo de Investigación de Análisis Matemático y Aplicaciones (AMA), Departamento de Matemática Aplicada a la Ingeniería Industrial, Escuela Técnica Superior de Ingeniería y Diseño Industrial, Universidad Politécnica de Madrid, Ronda de Valencia 3, 28012 Madrid, Spain
    \newline
		\href{https://orcid.org/0000-0002-1016-096X}{ORCID: \texttt{0000-0002-1016-096X} }}
	\email{\texttt{dl.rodriguez.vidanes@upm.es}}


\begin{document}

\begin{abstract}
We study norm-attaining operators on three renormings of $c_0$: Read's space $\Read$, its smooth renorming $\sRead$, and Mart\'in's space $\X$. We prove that, for every Banach space $F$ and every $Z\in\{\Read,\sRead\}$, norm-attaining finite-rank operators are dense in $\K(Z,F)$. Consequently, $\NA(Z,F)$ is dense in $\Lin(Z,F)$ whenever $F$ is reflexive. More generally, we establish that the compact Hilbertian smoothing construction of Debs, Godefroy and Saint-Raymond \cite{DGS} preserves the density of norm-attaining operators for all finite-dimensional ranges. We also extend a second-order criterion of Kadets, L\'opez-P\'erez, Mart\'in and Werner \cite{KLMWfr} from rank two to arbitrary ranks $k\geq 2$ for operators into $\ell_2^k$, and prove that the same smoothing preserves the existence and non-existence of such operators. We further investigate the possibilities and limitations of renorming techniques for excluding rank-two norm-attaining operators through analytical and geometrical approaches, including constructions that preserve the set of norm-attaining functionals. For Mart\'in's space $\X$, we prove that, for every Banach space $F$ and every integer $1\leq k\leq\dim F$, there exists a rank $k$ norm-attaining operator from $\X$ into $F$, answering a question of Mart\'in \cite{Mar}. We also show that, for every nonzero Banach space $F$, $\K(\X,F)$ contains a norm-dense linear subspace whose nonzero elements fail to attain their norms.
\end{abstract}

\maketitle

\section{Introduction}

All Banach spaces in this paper are real. For Banach spaces $Z$ and $F$, let $\Lin(Z,F)$, $\Fin(Z,F)$, and $\K(Z,F)$ denote, respectively, the spaces
of bounded, finite-rank, and compact linear operators from $Z$ to $F$, all endowed with the operator norm. We write $Z^*:=\Lin(Z,\R)$. An
operator $T\in\Lin(Z,F)$ is said to \emph{attain its norm} if $\|Tx\|=\|T\|$ for some unit vector $x\in Z$. The set of all such operators is denoted by $\NA(Z,F)$.

The Bishop--Phelps theorem \cite{BP} states that $\NA(Z,\R)$ is norm dense in $Z^*$ for every Banach space $Z$. Lindenstrauss \cite{Lin} showed that the analogous assertion for operators fails in general and introduced \emph{property B}: a Banach space $F$ has property B if $\NA(Z,F)$ is dense in $\Lin(Z,F)$ for every Banach space $Z$. He also proved that property $\beta$ implies property B \cite[Proposition~3]{Lin}. All known counterexamples to property B involve infinite-dimensional ranges. In particular, Gowers \cite{Gow} proved that $\ell_p$ fails property B for $1<p<\infty$, Acosta \cite{Aco} extended this conclusion to every infinite-dimensional strictly convex Banach space, and Mart\'in \cite{MarCompact} constructed compact operators that cannot be approximated by norm-attaining operators. By contrast, whether every finite-dimensional Banach space has property B remains open, even for the Euclidean plane,
a case explicitly raised by Johnson and Wolfe \cite[Question~6]{JW}. Since every operator into a finite-dimensional space is compact, property B and its compact-operator counterpart coincide for such ranges.

Motivated by the latter question, we study norm-attaining operators on three renormings of the classical sequence space $c_0$: Read's space $\Read$ \cite{Read}, its smooth renorming $\sRead$ \cite{KLMRead,DGS}, and Mart\'in's space $\X$ \cite{Mar}. Our main aim is to determine how the restricted linear structure of their norm-attaining functionals affects the existence and density of
norm-attaining operators of higher rank, and whether these spaces can serve as domains for counterexamples to property B with finite-dimensional ranges.

Read \cite{Read} constructed $\Read$ so that it has no proximinal subspaces of codimension two, and Rmoutil \cite{Rmo} subsequently proved
that $\NA(\Read,\R)$ contains no two-dimensional linear subspace. Using the smoothing construction of Debs, Godefroy and Saint-Raymond \cite{DGS}, Kadets, L\'opez-P\'erez and Mart\'in obtained $\sRead$,
a smooth renorming that preserves the set of norm-attaining functionals \cite[Lemma~11 and Example~12]{KLMRead}. Nevertheless, Kadets, L\'opez-P\'erez, Mart\'in and Werner \cite[Corollary~13.10 and Example~13.13]{KLMWfr} proved that both $\Read$ and $\sRead$ admit rank-two norm-attaining operators into every Banach space of dimension at least two.

Mart\'in \cite{Mar} constructed another smooth renorming $\X$ of $c_0$ with an even more restrictive structure of norm-attaining functionals: for every two-dimensional subspace $E\subset\X^*$,
the set $\NA(\X,\R)\cap E$ is contained in the union of two one-dimensional subspaces, and hence contains no cone generated by two linearly independent functionals. This suggested a possible approach to the property B problem for the Euclidean plane. Indeed, if every norm-attaining operator from $\X$ into $\ell_2^2$ had rank at most one, then $\ell_2^2$ would fail property B, since the operators of rank at most one form a proper norm-closed subset of $\Lin(\X,\ell_2^2)$. Mart\'in therefore asked whether there exists a rank-two norm-attaining operator from $\X$ into $\ell_2^2$ \cite[Open Problem~7]{Mar}.

The existence of linear subspaces within sets of norm-attaining and non-norm-attaining functionals has been extensively investigated, including questions concerning their dimension, density, and closedness \cite{AAAGP,BG,FI,GP,GPP,Mar}. Related questions have also been studied for norm-attaining and non-norm-attaining operators \cite{DFJRV,FGMR,PT}. A useful sufficient condition helps put these restrictions in perspective. We include an elementary proof, in Lemma~\ref{lem:densesubspace}, of a consequence of \cite[Corollaries~13.6 and 13.25]{KLMWfr}: if $\NA(Z,\R)$ contains a norm-dense linear subspace of $Z^*$, then norm-attaining finite-rank operators from $Z$ into any Banach space $F$ are dense among finite-rank operators. If, in addition, $Z^*$ has the approximation property, they are also dense in $\K(Z,F)$, since finite-rank operators are dense in $\K(Z,F)$ under the same assumption (see \cite[Theorem~1.e.5]{LT}). The duals of all three renormings considered here are isomorphic to $\ell_1$, and therefore have the approximation property. However, the restrictions described above rule out the existence of a norm-dense linear subspace of norm-attaining functionals. Thus, this sufficient condition does not apply to any of our three spaces, although its failure does not imply failure of density.

Our first tool is a rank-preserving construction for additive renormings of the form $p=q+s$, where $q$ is a norm on a Banach space and $s$ is a $q$-continuous seminorm (Lemma~\ref{lem:lift}). A norm-attaining operator for $q$ is converted into one for $p$ by a rank-one perturbation that preserves its range; when $s$ is small relative to $q$, the perturbation is small as well. This leads to a density-transfer principle for such renormings and,
when applied to finite truncations of Read's norm, yields
    $$
    \overline{\NA(\Read,F)\cap\Fin(\Read,F)}=\K(\Read,F)
    $$
for every Banach space $F$ (Theorem~\ref{thm:read-finite-density}). In particular, the restricted structure of $\NA(\Read,\R)$ does not prevent norm-attaining finite-rank operators from being dense among
compact operators.

We next establish a different transfer principle for smooth renormings. We prove that if norm-attaining operators from a Banach space are dense for every finite-dimensional range, then this property is preserved by the compact Hilbertian smoothing construction \cite{DGS} of Debs, Godefroy and Saint-Raymond
(Theorem~\ref{thm:transfer}). Combined with the result for $\Read$ and finite-rank approximation of compact operators, this gives
    $$
    \overline{\NA(\sRead,F)\cap\Fin(\sRead,F)}=\K(\sRead,F)
    $$
for every Banach space $F$ (Corollary~\ref{cor:smoothread}). Since $\Read$ and $\sRead$ are isomorphic to $c_0$, every bounded linear operator from either space into a reflexive Banach space is compact (see \cite[Theorem~2.4.10]{AK}). Consequently, $\NA(Z,F)$ is dense in $\Lin(Z,F)$ for $Z\in\{\Read,\sRead\}$ whenever $F$ is reflexive. In particular, neither space can serve as a domain for a counterexample to property B with a finite-dimensional range.

We also investigate the existence of norm-attaining operators of prescribed rank into Euclidean spaces. Extending the rank-two criterion of Kadets, L\'opez-P\'erez, Mart\'in and Werner \cite[Theorem~13.46 and Corollary~13.47]{KLMWfr} to every integer $k\geq 2$, we first obtain a characterization in terms of norm-attaining functionals and suitable independent directions (Proposition~\ref{prop:quadratic-characterization}). We then express this condition through the second-order behavior of the dual norm, proving that rank $k$ norm-attaining operators into $\ell_2^k$ exist precisely when the space of quadratic directions at some norm-attaining functional has dimension at least $k$ (Propositions~\ref{prop:theta-dual} and~\ref{prop:quadratic-directions} and Corollary~\ref{cor:dual-quadratic-nonexistence}). Using this characterization, we show that the compact Hilbertian smoothing construction preserves both the existence and the non-existence of such operators for every $k$ (Corollary~\ref{cor:dgs-rank-k}). Finally, we construct compact enlargements of the unit ball that preserve the set of norm-attaining functionals while eliminating all nontrivial quadratic directions at a prescribed norm-attaining ray (Proposition~\ref{prop:one-ray-elimination}).

For Mart\'in's space $\X$, the rank-preserving construction for additive renormings provides a positive answer to the question in \cite{Mar}. In fact, we prove that for every Banach space $F$ and every positive integer $k\leq\dim F$, where $\dim$ denotes the algebraic dimension, there exists an operator of rank $k$ in $\NA(\X,F)$ (Theorem~\ref{thm:finite-short}). This answers \cite[Open Problem~7]{Mar} affirmatively and \cite[Problem~13.12]{KLMWfr} negatively. We also establish a related rank-two preservation result for suitable absolute combinations of a norm and a continuous seminorm (Proposition~\ref{prop:absolute-lifting}), including $\ell_r$-type combinations (Corollary~\ref{cor:ell-r-lifting}).

We further investigate the limitations of Mart\'in's construction as a possible approach to excluding rank-two norm-attaining operators. We first show that whenever $\NA(Z,\R)$ contains a norm-dense linear subspace of $Z^*$ of dimension at least two, the norm-attaining functionals that arise as first coordinates of rank-two norm-attaining operators into $\ell_2^2$ form a norm-dense subset of the dual unit sphere (Proposition~\ref{prop:density-bad-coordinates}). In particular, this applies to the base space used in Mart\'in's construction (Corollary~\ref{cor:dense-bad-coordinates}). We then prove that its unit ball cannot be properly enlarged by a closed absolutely convex set while preserving the norms of all these functionals (Corollary~\ref{cor:no-simultaneous-repair-martin}). Together with our rank-preserving construction, these results show that modifying the auxiliary components of Mart\'in's construction while retaining its base norm and additive form cannot eliminate all rank-two norm-attaining operators. Turning to $\X$ itself, we prove that, for every nonzero Banach space $F$, there exists a norm-dense linear subspace of $\K(\X,F)$ whose nonzero elements fail to attain their norms (Proposition~\ref{prop:lineable}). However, this does not imply that norm-attaining operators are not dense. 

The paper is organized accordingly: Section~\ref{sec:prelim} introduces the notation and the three domain spaces, Section~\ref{sec:additive} treats additive renormings and Read's space, Section~\ref{sec:smooth} treats compact Hilbertian smoothing, Section~\ref{sec:rank-k} develops the Euclidean rank criteria, and Section~\ref{sec:martin} studies Mart\'in's space and its remaining density problem.

\section{Preliminaries}\label{sec:prelim}

In this section, we introduce the notation, recall the main features of the three renormings of $c_0$ considered in this paper, and establish auxiliary results needed in the subsequent sections.

\subsection{Notation}
For a Banach space $Z$, we write $B_Z$ and $S_Z$ for the closed unit ball and unit sphere, respectively. 
When a space carries several equivalent norms, we index the relevant objects by the norm: for a norm $r$ on $Z$, $B_r$, $S_r$ and $r^*$ denote the closed unit ball, the unit sphere and the dual norm of $(Z,r)$, and for $T\in\Lin(Z,F)$, we write $\|T\|_r$ for the operator norm of $T$ with respect to $r$ on $Z$, and $\|T\|_{r\to\rho}$ when the norm $\rho$ on the range has to be indicated as well. Observe that for equivalent norms $p,q$, $\|T\|_q\le C\|T\|_p$ whenever $p\le Cq$. For $T\in\Lin(Z,F)$, we write $\ran T:=T(Z)$
for its range, $\rank T:=\dim(\ran T)$
for its rank, and we put
    $$
    d_r(T):=\inf\{\|T-A\|_r \colon A\in\NA((Z,r),F)\}
    $$
and
    $$
    d_r^\Fin(T):=\inf\{\|T-A\|_r \colon A\in\NA((Z,r),F)\cap\Fin(Z,F)\}.$$
Note that both are finite since the zero operator belongs to $\NA((Z,r),F)\cap\Fin(Z,F)$. Throughout this paper, $P_M$ is the canonical projection from $c_0$ onto the first $M$ coordinates, and $E_M:=P_M(c_0)$.
Finally, recall that a norm $q$ on a Banach space $Z$ is \emph{smooth} if it is G\^ateaux differentiable at every non-zero point, equivalently, if every $x\ne0$ has exactly one support functional.

\subsection{Read's space}\label{subsec:readspace}

We recall Read's construction from \cite[Section~1]{Read}, using the normalized representation in \cite[Section~1]{KLMRead}. 
\begin{enumerate}
\item[(R1)] Let $c_{00}(\mathbb Q)$ denote the finitely supported sequences with rational entries. Let $(u_n)_{n}$ be a sequence in $c_{00}(\mathbb Q)$ which lists every element infinitely many times, and let $(a_n)_n$ be a strictly increasing sequence of positive integers such that
    $$
    a_n>\max(\supp u_n)\quad\text{and}\quad a_n\ge\|u_n\|_1,
    $$
where, for any $x=(x_n)_n\in c_{00}(\mathbb Q)$, we denote $\supp x := \{n\in \mathbb N \colon x_n\neq 0 \}$.

\item[(R2)] Regard each $u_n$ as an element of $\ell_1$ (the dual of $c_0$), and let $(e_j)_j$ be the canonical basis for $\ell_1$. Define
    $$
    p_{\Read}(x):=\|x\|_\infty+\sum_{n=1}^{\infty}2^{-a_n^2}|(u_n-e_{a_n})(x)|,
    $$
for every $x\in c_0$.
Read's space is $\Read:=(c_0,p_{\Read})$.

\item[(R3)] Set
    $$
    v_n:=\frac{u_n-e_{a_n}}{\|u_n-e_{a_n}\|_1},\quad r_n:=2^{-a_n^2}\|u_n-e_{a_n}\|_1.
    $$
Then $v_n\in c_{00}$, $\|v_n\|_1=1$, and $\supp v_n\subset\{1,\ldots,a_n\}$. Also, we see that
\begin{equation}\label{eq:read-norm}
p_{\Read}(x)=\|x\|_\infty+\sum_{n=1}^{\infty}r_n|v_n(x)|,
\end{equation}
for every $x\in c_0$.

\item[(R4)] It holds that $\sum_{n=1}^\infty r_n\le2$. Therefore, $\|x\|_\infty\le p_{\Read}(x)\le3\|x\|_\infty$ for every $x\in c_0$. Hence, $p_{\Read}$ is an equivalent norm on $c_0$.
\end{enumerate}

For the density argument in Section~\ref{sec:additive}, define, for $N\in\N$ and $x\in c_0$,
    $$
    q_N(x):=\|x\|_\infty+\sum_{n\le N}r_n|v_n(x)|,\quad s_N(x):=\sum_{n>N}r_n|v_n(x)|,\quad\text{and}\quad\eta_N:=\sum_{n>N}r_n.
    $$

\begin{lemma}\label{lem:trunc}
    For every $N\in\N$, the following hold.
        \begin{enumerate}
        \item[\rm{(a)}] $p_\Read=q_N+s_N$, $s_N\le\eta_Nq_N$, and $\eta_N\to0$ as $N\to\infty$.
        \item[\rm{(b)}] $P_M$ is $q_N$-contractive for every integer $M\ge a_N$.
        \end{enumerate}
\end{lemma}

\begin{proof}
Since $\|v_n\|_1=1$ and $q_N\ge\|\cdot\|_\infty$, we have $s_N\le\eta_N\|\cdot\|_\infty\le\eta_Nq_N$. The remaining assertions in (a) follow from the definitions and the summability of $(r_n)_n$.

For (b), if $M\ge a_N$ and $n\le N$, then $v_n(P_Mx)=v_n(x)$ because $\supp v_n\subseteq\{1,\ldots,a_N\}$. Together with $\|P_Mx\|_\infty\le\|x\|_\infty$, this gives $q_N(P_Mx)\le q_N(x)$.
\end{proof}

\subsection{Smooth Read spaces}\label{subsec:smoothread}

We use the smoothing construction from the proof of \cite[Theorem~9(4)]{DGS}. The argument in \cite[Lemma~11]{KLMRead} gives the following formulation.

\begin{lemma}\label{lem:smoothing}
    Let $(X,q)$ be a Banach space, let $G$ be a Hilbert space, and let $U:G\to(X,q)$ be compact with dense range. Then,
        \begin{equation}\label{eq:smoothing-ball}
        B_p:=B_q+U(B_G)
        \end{equation}
    defines an equivalent smooth norm $p$ on $X$. Its dual norm is
        \begin{equation}\label{eq:smoothing-dual}
        p^*(f)=q^*(f)+\|U^*f\|_G,
        \end{equation}
    for every $f\in X^*$, and $\NA((X,p),\R)=\NA((X,q),\R)$.
\end{lemma}

Let $G$ be a Hilbert space and fix a compact operator $U:G\to\Read$ with dense range. We write $\sRead$ for the space with unit ball $B_{p_\Read}+U(B_G)$.

By Lemma~\ref{lem:smoothing}, the norm-attaining functionals of $\sRead$ are exactly those of $\Read$. In particular, they contain no two-dimensional linear subspace \cite[Theorem~4.2]{Rmo}.

\subsection{Mart\'in's space \texorpdfstring{$\X$}{X}}\label{subsec:space}
We recall the construction in \cite[Theorem~1]{Mar}. The items below summarize Steps~1--7 of its proof. Let $(e_{n,m})_{n,m}$ be the canonical basis of $\ell_1(\N\times\N)$.
\begin{enumerate}
\item[(M1)] $c_0$ is endowed with an equivalent smooth norm $\tn{\cdot}$ such that $$\NA((c_0,\tn{\cdot}),\R)=\NA(c_0,\R),$$
and $Y$ is a separable operator range, norm dense in $(c_0,\tn{\cdot})^*$, with 
    $$
    \NA((c_0,\tn{\cdot}),\R)\cap Y=\{0\}.
    $$
\item[(M2)] $T:\ell_1(\N\times\N)\to Y\subseteq (c_0,\tn{\cdot})^*$ is a norm-one injective operator such that, writing $v^*_{n,m}:=T(e_{n,m})$, for every $m\in\N$ the set $\{v^*_{n,m}/\tn{v^*_{n,m}}^*:\ n\in\N\}$ is norm dense in $S_Y$, hence in $S_{(c_0,\tn{\cdot})^*}$ \cite[Lemma~B]{Mar}.
\item[(M3)] For $m,n\in\N$, put $\Phi_m(n):=2^{-m}2^{-n}\tn{v^*_{n,m}}^*$. Let $|\cdot|_m$ be the equivalent norm on $\ell_1$ given by \cite[Lemma~A]{Mar} for the sequence $(\Phi_m)_m$.
\item[(M4)] Define the operator $R_m:(c_0,\tn{\cdot})\to(\ell_1,|\cdot|_m)$ by
$$[R_m(x)](n)=\frac{m}{2^m}\frac{1}{2^n}v^*_{n,m}(x),$$
for every $n\in\N$ and $x\in c_0$.
Then $R_m$ is injective and $|R_m(x)|_m\le\|R_m(x)\|_1\le \frac{m}{2^m}\tn{x}$.
\item[(M5)] Set $W:=\big(\bigoplus_{m\in\N}(\ell_1,|\cdot|_m)\big)_{\ell_1}$ and the operator $R:(c_0,\tn{\cdot})\to W$, by $R(x):=(R_m(x))_{m}$.
Then $\|Rx\|_W=\sum_{m=1}^\infty|R_m x|_m\le 2\tn{x}$ and $R$ is injective. Also $W^*=\big(\bigoplus_{m\in \mathbb N}(\ell_\infty,|\cdot|_m^*)\big)_{\ell_\infty}$, so that $\|(w_m)_m\|_{W^*}=\sup_m|w_m|_m^*$.
\item[(M6)] We have $R^*(W^*)\subseteq Y$. More precisely, for $w=(w_m)_m\in W^*$,
$$R^*w=T\left(\sum_{m,n}w_m(n)\frac{m}{2^{m+n}}e_{n,m}\right)\in Y.$$
The series inside $T$ converges absolutely in $\ell_1(\N\times\N)$, since
$|w_m(n)|\le |w_m|_m^*\le\|w\|_{W^*}$.
\item[(M7)] Put $p(x):=\tn{x}+\|Rx\|_W$ for $x\in c_0$, and $\X:=(c_0,p)$. Then $\tn{x}\le p(x)\le 3\tn{x}$, and $\X$ is smooth. Furthermore, if $x\in S_\X$ and $f\in S_{\X^*}$ satisfy $f(x)=1$, then $f=f_0+R^*w$ for some $f_0\in S_{(c_0,\tn{\cdot})^*}$ and $w=(w_m)_m\in W^*$ such that
    $$
    f_0(x)=\tn{x},\quad w(Rx)=\|Rx\|_W,\quad|w_m|_m^*=1\quad\text{and}\quad w_m(R_mx)=|R_mx|_m,
    $$
for every $m\in\N$.
\end{enumerate}

\section{Read's space}\label{sec:additive}

In this section, we establish general results on norm-attaining operators under additive renormings, together with auxiliary density and approximation results. We then apply them to Read's space, proving that norm-attaining finite-rank operators are dense among compact operators into every Banach space and, consequently, that norm-attaining operators are dense among all bounded operators whenever the range is reflexive.

Throughout this section, $(Z,q)$ and $F$ are Banach spaces, $s$ is a $q$-continuous seminorm on $Z$, that is, $s\le Cq$ for some $C\ge0$, and set $p:=q+s$. Then $p$ is a norm on $Z$ equivalent to $q$, with $q\le p\le(1+C)q$.
Also, an operator $P\in\Lin((Z,q),(Z,q))$ is said to be \textit{$q$-contractive}
if $q(Px)\leq q(x)$ for every $x\in Z$.

\begin{lemma}\label{lem:lift}
Let $S\in\Lin((Z,q),F)\setminus\{0\}$ attain its norm at $z_0\in S_q$, and let $h\in Z^*$ satisfy $|h|\le s$ and $h(z_0)=s(z_0)$ (such a $h$ exists by the Hahn--Banach theorem). Then
    $$
    A:=S+h\otimes Sz_0 \in \Lin((Z,q),F)
    $$
attains its norm at $\tfrac{z_0}{p(z_0)}\in S_p$, $\|A\|_p=\|S\|_q$, $\ran A=\ran S$ and $\rank(A-S)\le1$. If, moreover, $s\le\eta q$ for some $\eta\ge0$, then $\|A-S\|_q\le\eta\|S\|_q$.
\end{lemma}
\begin{proof}
For $z\in Z$, we have
    $$
    \|Az\|\le\|Sz\|+|h(z)|\|Sz_0\|\le\|S\|_q(q(z)+s(z))=\|S\|_qp(z).
    $$
Hence, by taking the supremum over $z\in B_p$, we obtain
$$\|A\|_p\le\|S\|_q.$$ 
Also, from $Az_0=(1+s(z_0))Sz_0=p(z_0)Sz_0$, we have
    $$
    \left\|A\left(\tfrac{z_0}{p(z_0)}\right)\right\|=\|Sz_0\|=\|S\|_q.
    $$

Moreover, $A=S\circ(\id{Z}+h\otimes z_0)$. Since $h(z_0)=s(z_0)\ge0$, the operator $\id{Z}+h\otimes z_0$ is invertible with inverse $\id{Z}-(1+h(z_0))^{-1}h\otimes z_0$. Hence, $\ran A=\ran S$ and $A-S=h\otimes Sz_0$ has rank at most one. 

Finally, if $s\le\eta q$, then for every $z\in Z$,
    $$
    \|(A-S)z\|=|h(z)|\|Sz_0\|\le s(z)\|S\|_q\le\eta q(z) \|S\|_q.
    $$
Taking the supremum over $z\in B_q$ gives the desired estimate.
\end{proof}

\begin{theorem}\label{thm:permanence}
Let $(Z,p)$ and $F$ be Banach spaces. Assume that, for each $\eta>0$, there exist a norm $q_\eta$ and a seminorm $s_\eta$ on $Z$ such that 
$$p=q_\eta+s_\eta\quad\text{and}\quad s_\eta\le\eta q_\eta.$$
Then, the following assertions hold.
\begin{enumerate}
    \item[\rm{(a)}] If $\lim_{\eta\to0}d_{q_\eta}(T)=0$, for every $T\in\Lin((Z,p),F)$, then $\NA((Z,p),F)$ is dense in $\Lin((Z,p),F)$. 
    \item[\rm{(b)}] If $\lim_{\eta\to0}d_{q_\eta}^\Fin(T)=0$, for every $T\in\K((Z,p),F)$, then $\NA((Z,p),F)\cap\Fin((Z,p),F)$ is dense in $\K((Z,p),F)$. 
\end{enumerate}
\end{theorem}
\begin{proof}
For every $\eta>0$, the assumptions give
$$q_\eta\le p= q_\eta+s_\eta\le(1+\eta)q_\eta.$$
Thus $q_\eta$ is an equivalent norm on $Z$.

For any $S\in\NA((Z,q_\eta),F)\setminus\{0\}$, Lemma~\ref{lem:lift} gives $A\in\NA((Z,p),F)$ such that
$\|A-S\|_p\le\eta\|S\|_{q_\eta}$. For $S=0$, take $A=0$. Therefore,
    \begin{align*}
        d_p(T)&\le\|A-T\|_p\\
        &\le\|A-S\|_p+\|S-T\|_p \\
        &\le\eta\|S\|_{q_\eta}+\|S-T\|_{q_\eta} \\
        &\le\eta\|T\|_{q_\eta}+(1+\eta)\|S-T\|_{q_\eta}.
    \end{align*}
Taking the infimum over $S\in\NA((Z,{q_\eta}),F)$, and using $\|T\|_{q_\eta}\le(1+\eta)\|T\|_p$, we obtain
$$d_p(T)\le\eta(1+\eta)\|T\|_p+(1+\eta)d_{q_\eta}(T).$$
Letting $\eta\to0$ gives $d_p(T)=0$.

Since the construction in Lemma~\ref{lem:lift} preserves finite rank, (b) follows by the same argument.
\end{proof}

\begin{lemma}\label{lem:projections}
Let $q$ be an equivalent norm on $c_0$ and $F$ be a Banach space. For every $T\in\K(c_0,F)$,
\begin{enumerate}
\item[\rm{(a)}] $\|T-TP_M\|_q\to0$ as $M\to\infty$.
\item[\rm{(b)}] If $P_M$ is $q$-contractive, then $TP_M\in\NA((c_0,q),F)\cap\Fin(c_0,F)$.
\end{enumerate}
Consequently, if $P_M$ is $q$-contractive for arbitrarily large integers $M$, then
$$\overline{\NA((c_0,q),F)\cap\Fin((c_0,q),F)}=\K((c_0,q),F).$$
\end{lemma}

\begin{proof}
For (a), first endow $c_0$ with its usual norm. Since $T$ is compact, $T^*(B_{F^*})$ is relatively compact in $\ell_1$ by Schauder's theorem (see, for instance, \cite[Theorem~3.4, p. 174]{C}). The coordinate tails converge uniformly on relatively compact subsets of $\ell_1$, so
$$\|T-TP_M\|_\infty=\sup_{f\in B_{F^*}}\|(\id{\ell_1}-P_M^*)T^*f\|_1\longrightarrow0.$$
The equivalence of $q$ and $\|\cdot\|_\infty$ now gives (a).

To show (b), put $S:=TP_M$.
For $x\in c_0$,
    $$
    \|Sx\|=\|S(P_Mx)\|\le\|S|_{E_M}\|_qq(P_Mx)\le\|S|_{E_M}\|_qq(x).
    $$
Hence $\|S\|_q=\|S|_{E_M}\|_q$. Since the unit sphere of $(E_M,q)$ is compact, there is $z_0\in E_M$ with $q(z_0)=1$ and $\|Sz_0\|=\|S|_{E_M}\|_q=\|S\|_q$, hence $S\in\NA((c_0,q),F)$.

The last assertion follows from (a) and (b).
\end{proof}

The following consequence of \cite[Corollaries~13.6 and~13.25]{KLMWfr} is included with its elementary finite-rank proof.

\begin{lemma}\label{lem:densesubspace}
Let $Z$ and $F$ be Banach spaces. Assume that $D$ is a linear subspace of $Z^*$ with $D\subseteq\NA(Z,\R)$. Then the following assertions hold:
\begin{itemize}
    \item[(a)] Every finite-rank operator $S:Z\to F$ with $S^*(F^*)\subseteq D$ attains its norm.
    \item[(b)] If $D$ is dense in $Z^*$, then the operators in \textup{(a)} are dense in $\Fin(Z,F)$.
\end{itemize}
\end{lemma}
\begin{proof}
For (a), assume $S\ne0$ and regard $S$ as taking values in the finite-dimensional space $\ran S$.
The adjoint $S^*:(\ran S)^*\to Z^*$ attains its norm at some $u^*\in S_{(\ran S)^*}$.
By Hahn--Banach theorem, extending $u^*$ to $F^*$ shows that $S^*u^*\in D$.
Thus, there exists $z\in S_Z$ such that 
    $$
    \|S\| = \|S^*u^*\| = |(S^*u^*)(z)| = |u^*(Sz)| \le \|Sz\| \le \|S\|.
    $$
In particular, $\|Sz\|=\|S\|$.

To show (b), let $S:Z\to F$ be a nonzero finite-rank operator and $\varepsilon>0$.
Choose a normalized basis $u_1,\ldots,u_k$ of $\ran S$, with coordinate functionals $u_1^*,\ldots,u_k^*$ in $(\ran S)^*$. Regarding $S$ as an operator into $\ran S$, put $h_j=u_j^*\circ S\in Z^*$. Then $S=\sum_{j=1}^k h_j\otimes u_j$.
Since $D$ is dense in $Z^*$, take $g_j\in D$ such that $\|g_j-h_j\|<\tfrac\varepsilon k$ for each $1\le j\le k$, and define $T=\sum_{j=1}^k g_j\otimes u_j$. Then we have
$$\|T-S\|\le\sum_{j=1}^k\|g_j-h_j\|\|u_j\|<\varepsilon.$$
Moreover, for every $f\in F^*$, $T^*f=\sum_{j=1}^k f(u_j)g_j\in D$. Thus $T$ attains its norm by the first assertion.
\end{proof}

Finite truncations of Read's norm and Lemma~\ref{lem:lift} yield the following density result.

\begin{theorem}\label{thm:read-finite-density}
Let $\Read$ be Read's space as in Subsection~\ref{subsec:readspace}. For every Banach space $F$,
\begin{equation}\label{eq:readna}
    \overline{\NA(\Read,F)\cap\Fin(\Read,F)}=\K(\Read,F).
\end{equation}
In particular, $\NA(\Read,F)$ is dense in $\Lin(\Read,F)$ whenever $F$ is reflexive.
\end{theorem}
\begin{proof}
The inclusion $\subseteq$ is trivial. By (R4) and Lemma~\ref{lem:trunc}(b), $\|\cdot\|_\infty\le q_N\le3\|\cdot\|_\infty$ and $P_M$ is $q_N$-contractive for every $M\ge a_N$. Hence Lemma~\ref{lem:projections} gives
$$\overline{\NA((c_0,q_N),F)\cap\Fin((c_0,q_N),F)}=\K((c_0,q_N),F).$$
For each $\eta>0$, choose $N$ with $\eta_N\le\eta$ and take $q_\eta=q_N$ and $s_\eta=s_N$. Lemma~\ref{lem:trunc}(a) and Theorem~\ref{thm:permanence}(b) yield \eqref{eq:readna}.

If $F$ is reflexive and $T\in \Lin(\Read,F)$, then $T^*$ is weakly compact. Since $\Read^*$ is isomorphic to $\ell_1$, it has the Schur property, and hence $T^*$ is compact. Thus $T$ is compact by Schauder's theorem, which proves the final assertion.
\end{proof}

Theorem~\ref{thm:read-finite-density} excludes Read's space as a domain for a counterexample to Lindenstrauss property B of a $k$-dimensional range for $k\geq 2$.

\section{Smooth renormings}\label{sec:smooth}

In this section, we prove that density of norm-attaining operators for all finite-dimensional ranges is preserved by compact Hilbert smoothings. 

Let $(X,q)$ be a Banach space, let $G$ denote a Hilbert space with inner product $\ip{\cdot}{\cdot}$, and let $U:G\to(X,q)$ be compact with dense range. Throughout this section, $p$ is the smooth norm given by $B_p:=B_q+U(B_G)$ in Lemma~\ref{lem:smoothing}, and we identify $G^*$ with $G$.
 
We first record two geometric lemmas. For a Banach space $E$, a point $a\in S_E$ and $\lambda>0$, the set $w+\lambda(B_E-a)$ is the ball of radius $\lambda$ centered at $w-\lambda a$; its boundary passes through $w$. 
The first lemma places the image of a Hilbert unit ball inside such a ball, with $w$ the point of the image at which a prescribed functional is maximal. The second lemma turns an inclusion of this kind into a norm-attaining rank-one perturbation, and it holds for arbitrary Banach ranges.
 
\begin{lemma}\label{lem:nearball}
Let $G$ be a Hilbert space with inner product $\ip{\cdot}{\cdot}$, $L\in\Lin(G,G)$, and $\phi\in S_G$ with $L^*\phi\ne0$. Put $\psi:=L^*\phi/\|L^*\phi\|_G$. If
$$\lambda\ge\frac{\|L\|^2}{\|L^*\phi\|_G},$$
then
\begin{equation}\label{eq:nearball}
L(B_G)\subseteq L\psi+\lambda(B_G-\phi).
\end{equation}
\end{lemma}
\begin{proof}
For $g\in B_G$, we have 
    $$
    \|g-\psi\|_G^2=\|g\|_G^2+1-2\ip{\psi}{g}_G\le2(1-\ip{\psi}{g}_G)
    $$ 
and
    $$
    \ip{L(g-\psi)}{\phi}_G=\ip{g-\psi}{L^*\phi}_G=\|L^*\phi\|_G\left(\ip{\psi}{g}_G-1\right).
    $$ 
Hence,
    \begin{align*}
    \|L(g-\psi)\|_G^2
    &\le2\|L\|^2(1-\ip{\psi}{g}_G)
    \le2\lambda\|L^*\phi\|_G(1-\ip{\psi}{g}_G)
    =-2\lambda\ip{L(g-\psi)}{\phi}_G .
    \end{align*}
This yields $\|Lg-L\psi+\lambda\phi\|_G\le\lambda$, which is \eqref{eq:nearball}.
\end{proof}
 
\begin{lemma}\label{lem:recenter}
Let $Z$ and $E$ be Banach spaces, $S\in\Lin(Z,E)$, $z_0\in B_Z$, $a\in S_E$, $\xi\in S_{E^*}$ with $\xi(a)=1$, and $\lambda>0$. Assume that
\begin{equation}\label{eq:recenter-hyp}
S(B_Z)\subseteq Sz_0+\lambda(B_E-a)\quad\text{and}\quad \xi(Sz_0)>0.
\end{equation}
Then
    $$
    A=S+\frac{S^*\xi}{\xi(Sz_0)}\otimes(\lambda a-Sz_0) \in \Lin(Z,E)
    $$
satisfies $\|A\|=\lambda=\|Az_0\|$. In particular, $z_0\in S_Z$ and $A\in\NA(Z,E)$.
\end{lemma}
\begin{proof}
For every $z\in B_Z$, \eqref{eq:recenter-hyp} gives
$$\|Sz-Sz_0+\lambda a\|\le\lambda\quad\text{and}\quad\|Sz+Sz_0-\lambda a\|\le\lambda .$$
Applying $\xi$, we have $\xi(Sz)\le \xi(Sz_0)$ and $-\xi(Sz)\le\xi(Sz_0)$. Thus $|\xi(Sz)|\le \xi(Sz_0)$. Therefore $\|S^*\xi\|=\xi(Sz_0)$, so that $z_0\in S_Z$. For $x\in B_Z$, since $\tfrac{S^*\xi}{\xi(Sz_0)}(x)\in[-1,1]$,
\begin{align*}
    Ax&=Sx+\frac{S^*\xi}{\xi(Sz_0)}(x)(\lambda a-Sz_0)\\
    &=\frac{1+\tfrac{S^*\xi}{\xi(Sz_0)}(x)}{2}\left(Sx-Sz_0+\lambda a\right)+\frac{1-\tfrac{S^*\xi}{\xi(Sz_0)}(x)}{2}\left(Sx+Sz_0-\lambda a\right),
\end{align*}
hence $\|Ax\|\le\lambda$ by convexity. On the other hand, $Az_0=\lambda a$ shows that $\|A\|=\lambda=\|Az_0\|$.
\end{proof}
 
\begin{theorem}\label{thm:transfer}
Assume that $\NA((X,q),F)$ is dense in $\Lin((X,q),F)$ for every finite-dimensional Banach space $F$. Then, for the norm $p$ in \eqref{eq:smoothing-ball},
$$\overline{\NA((X,p),F)}=\Lin((X,p),F)$$
for every finite-dimensional Banach space $F$.
\end{theorem}
\begin{proof}
The zero operator already attains its norm. By homogeneity, it suffices to approximate $T\in\Lin((X,p),F)$ with $\|T\|_p=1$. Let $E=\ran T$ with the norm inherited from $F$. We regard $T$ as a surjection onto $E$ and construct the approximants with values in $E$; viewed as $F$-valued operators they have the same norms and the same norm-attaining points.

\medskip
 
\emph{Step 1:}
Since $T$ is continuous and maps $X$ onto $E$, and $U(G)$ is dense in $X$, the subspace $T(U(G))$ is dense in $E$. 
As $E$ is finite-dimensional, the map $TU$ is surjective, and hence $U^*T^*:E^*\to G$ is injective. Likewise, $T^*:E^*\to X^*$ is injective. By \eqref{eq:smoothing-dual} and the compactness of $S_{E^*}$,
\begin{equation}\label{eq:r}
1=\|T\|_p=\max_{f\in S_{E^*}}\left\{q^*(T^*f)+\|U^*T^*f\|_G\right\}>\max_{f\in S_{E^*}}\|U^*T^*f\|_G=\|TU\|,
\end{equation}
where the inequality is strict because $T^*$ is injective.
 
Define
$$D:=\{v\in E:\ v+TU(B_G)\subseteq B_E\}=\bigcap_{g\in B_G}\left(B_E-TUg\right).$$
Then $D$ is closed, convex and symmetric, and
\begin{equation}\label{eq:D-balls}
(1-\|TU\|)B_E\subseteq D\subseteq B_E\quad\text{and}\quad D+TU(B_G)\subseteq B_E.
\end{equation}
The first inclusion follows from the triangle inequality, the second from $0\in TU(B_G)$, and the last is the definition of $D$. Consequently, $D$ is the closed unit ball of a norm $p_D$ on $E$ equivalent to the given norm.
 
Since
$$T(B_q)+TU(B_G)=T(B_p)\subseteq B_E,$$
we have $T(B_q)\subseteq D$, that is, $\|T\|_{q\to p_D}\le1$. In fact, equality holds. Otherwise, $T(B_q)\subseteq tD$ for some $t\in[0,1)$. Writing $TUg=t\,TUg+(1-t)\,TUg$ for $g\in B_G$ and using \eqref{eq:D-balls}, we get
$$T(B_p)=T(B_q)+TU(B_G)\subseteq t\left(D+TU(B_G)\right)+(1-t)\,TU(B_G)\subseteq\left(t+(1-t)\|TU\|\right)B_E .$$
Hence $\|T\|_p\le t+(1-t)\|TU\|<1$ by \eqref{eq:r}, which contradicts $\|T\|_p=1$.

\medskip
 
\emph{Step 2:}
By the hypothesis applied to $(E,p_D)$, there is a sequence in $\NA((X,q),(E,p_D))$ converging to $T$. Since $\|T\|_{q\to p_D}=1$, normalizing its terms we obtain $(S_n)_n \subseteq \NA((X,q),(E,p_D))$ such that $\|S_n\|_{q\to p_D}=1$ and $\|S_n-T\|_{q\to p_D}\longrightarrow0$.
Since $p_D$ and the given norm of $E$ are equivalent, also $\|S_n-T\|_{q\to E}\to0$.

For every $n\in \mathbb N$, choose $x_n\in S_q$ with $p_D(S_nx_n)=1$, and put $v_n:=S_nx_n\in\partial D$. 
Since $S_nx_n\in D$, the compact set $S_nx_n+TU(B_G)$ is contained in $B_E$, and it meets $S_E$: otherwise $\max_{g\in B_G}\|S_nx_n+TUg\|_E<1$, and $D$ would contain a neighbourhood of $S_nx_n$, contradicting $S_nx_n\in\partial D$. Hence there are $\phi_n\in B_G$ and $\xi_n\in S_{E^*}$ such that
    \begin{equation}\label{eq:contact-point}
    a_n:=S_nx_n+TU\phi_n\in S_E\quad\text{and}\quad\xi_n(a_n)=1 .
    \end{equation}
Moreover, $\phi_n\in S_G$. Indeed, $\xi_n(S_nx_n+TUg)\le1=\xi_n(S_nx_n+TU\phi_n)$ for every $g\in B_G$, so the functional $\ip{U^*T^*\xi_n}{\cdot}_G$ attains its maximum over $B_G$ at $\phi_n$, and this functional is nonzero because $U^*T^*$ is injective.

\medskip
 
\emph{Step 3:}
We have
    \begin{equation}\label{equ:CnC}
        \|S_nU-TU\|\le\|S_n-T\|_{q\to E}\,\|U\|_{G\to q}\longrightarrow0.
    \end{equation}
Since $TU$ maps $G$ onto the finite-dimensional space $E$, there is a linear operator $J:E\to G$ such that $TUJ=\id{E}$. For each $n\in \mathbb N$, set
    $$
    L_n:=\id{G}+J(S_nU-TU).
    $$
Then
    $$
    TUL_n=S_nU\quad\text{and}\quad\|L_n-\id{G}\|\le\|J\|\,\|S_nU-TU\|\longrightarrow0.
    $$
For all sufficiently large $n$, we have $\|L_n-\id{G}\|<1$, and then $\|L_n^*\phi_n\|_G>0$. Define
    $$
    \psi_n:=\frac{L_n^*\phi_n}{\|L_n^*\phi_n\|_G}\in S_G\quad\text{and}\quad\lambda_n:=\max\left\{1,\frac{\|L_n\|^2}{\|L_n^*\phi_n\|_G}\right\}.
    $$
Since $\|L_n^*\phi_n\|_G\to1$ and $\|L_n\|\to1$, we have $\lambda_n\to1$, and
\begin{align}
    \|\psi_n-\phi_n\|_G&\le\|\psi_n-L_n^*\phi_n\|_G+\|L_n^*\phi_n-\phi_n\|_G\notag\\
    &\le\left|1-\|L_n^*\phi_n\|_G\right|+\|L_n-\id{G}\|\le2\|L_n-\id{G}\|\longrightarrow0 \label{equ:psinphin}.
\end{align}
Lemma~\ref{lem:nearball}, applied to $L_n$, $\phi_n$, $\psi_n$ and $\lambda_n$, gives $L_n(B_G)\subseteq L_n\psi_n+\lambda_n(B_G-\phi_n)$.
Applying $TU$ to both sides yields
\begin{equation}\label{eq:ellipsoid}
S_nU(B_G)\subseteq S_nU\psi_n+\lambda_n\left(TU(B_G)-TU\phi_n\right).
\end{equation}
Set
    $$
    z_n:=x_n+U\psi_n\in B_q+U(B_G)=B_p.
    $$
Then $S_nz_n=S_nx_n+S_nU\psi_n$.
Since $\lambda_n\ge1$ and $S_nx_n\in D$, convexity of $D$ gives
    \begin{equation}\label{eq:D-dilate}
    D\subseteq S_nx_n+\lambda_n(D-S_nx_n).
    \end{equation}
Indeed, for $v\in D$ the vector $\lambda_n^{-1}v+(1-\lambda_n^{-1})S_nx_n$ belongs to $D$. Since $S_n(B_q)\subseteq D$, adding \eqref{eq:D-dilate} and \eqref{eq:ellipsoid}, and using the definition of $a_n$ and \eqref{eq:D-balls}, we obtain
\begin{align}
S_n(B_p)=S_n(B_q)+S_nU(B_G)
&\subseteq S_nx_n+\lambda_n(D-S_nx_n)+S_nU\psi_n+\lambda_n\left(TU(B_G)-TU\phi_n\right)\notag\\
&=S_nz_n+\lambda_n\left(D+TU(B_G)-a_n\right)\subseteq S_nz_n+\lambda_n(B_E-a_n).\label{eq:tangent-ball}
\end{align}
This is the ball of radius $\lambda_n$ centered at $S_nz_n-\lambda_na_n$, and its boundary contains $S_nz_n$.

\medskip
 
\emph{Step 4:}
Since $S_nz_n-a_n=S_nU\psi_n-TU\phi_n=(S_nU-TU)\psi_n+TU(\psi_n-\phi_n)$, by \eqref{equ:CnC} and \eqref{equ:psinphin} we have
    \begin{equation}\label{equ:Snznan}
        \|S_nz_n-a_n\|\le\|S_nU-TU\|+\|TU\|\,\|\psi_n-\phi_n\|_G\longrightarrow0.
    \end{equation}
In particular, $\xi_n(S_nz_n)\ge\xi_n(a_n)-\|S_nz_n-a_n\|=1-\|S_nz_n-a_n\|>0$ for all sufficiently large $n$. For such $n$, by \eqref{eq:contact-point} and \eqref{eq:tangent-ball}, Lemma~\ref{lem:recenter} applies to $Z=(X,p)$, $S=S_n$, $z_0=z_n$, $a=a_n$, $\xi=\xi_n$ and $\lambda=\lambda_n$, and gives
    $$
    A_n:=S_n+\frac{S_n^*\xi_n}{\xi_n(S_nz_n)}\otimes(\lambda_na_n-S_nz_n) \in \NA((X,p),E)
    $$
with
    $$
    \|A_n-S_n\|_{p\to E}=\|\lambda_na_n-S_nz_n\|\le(\lambda_n-1)+\|a_n-S_nz_n\|\longrightarrow0
    $$
by \eqref{equ:Snznan} and $\lambda_n\to1$.
Since $p$ and $q$ are equivalent, $\|S_n-T\|_{p\to E}\longrightarrow0$ as well. Hence $\|A_n-T\|_{p\to E}\longrightarrow0$, and viewing $A_n$ as $F$-valued operators completes the proof.
\end{proof}
 
\begin{corollary}\label{cor:smoothread}
Let $\sRead$ be a smooth renorming of Read's space as in Subsection~\ref{subsec:smoothread}. For every Banach space $F$,
\begin{equation}\label{eq:KsRF}
    \overline{\NA(\sRead,F)\cap\Fin(\sRead,F)}=\K(\sRead,F).
\end{equation}
In particular, $\NA(\sRead,F)$ is dense in $\Lin(\sRead,F)$ whenever $F$ is reflexive.
\end{corollary}
\begin{proof}
Equation~\eqref{eq:KsRF} follows from
Theorems~\ref{thm:read-finite-density} and~\ref{thm:transfer},
together with Lemma~\ref{lem:projections}(a).

Since the dual of $\sRead$ is isomorphic to $\ell_1$, the same argument as in the proof of Theorem~\ref{thm:read-finite-density} shows that every bounded linear operator from $\sRead$ to a reflexive space is compact. This proves the final assertion.
\end{proof}

In particular, smooth Read spaces are not domains for counterexamples to Lindenstrauss property B with a $k$-dimensional range for $k\geq 2$.

\section{Norm-attaining operators onto $\ell_2^k$}\label{sec:rank-k}

We now turn from density to the existence of norm-attaining operators of a given rank into Euclidean spaces. 

For a normed space $X$ and an integer $k\geq2$, we write
    $$
    \NA^{(k)}(X,\ell_2^k):=\{T\in\NA(X,\ell_2^k) \colon \rank T=k\}.
    $$
Kadets, L\'opez-P\'erez, Mart\'in and Werner described the operators in $\NA^{(2)}(X,\ell_2^2)$ by means of mates and gave a second-order criterion for their existence \cite[Theorem~13.46 and Corollary~13.47]{KLMWfr}.

Let $m\geq1$ be an integer. For $f\in S_{X^*}$ and $g=(g_1,\ldots,g_m)\in(X^*)^m$, write $g(x):=(g_1(x),\ldots,g_m(x))\in\ell_2^m$ for $x\in X$, and define
    $$
    \vartheta(f,g):=\inf\left\{\frac{1-f(x)^2}{\|g(x)\|_2^2} \colon x\in B_X,\ g(x)\neq0\right\}.
    $$
We adopt the convention that the infimum of the empty set is $+\infty$. Since $|f(x)|\leq1$ for every $x\in B_X$, we have $\vartheta(f,g)\in[0,+\infty]$. For $m=1$, we identify $g$ with $g_1\in X^*$.

For $k=2$, the next proposition is a reformulation of \cite[Theorem~13.46 and Corollary~13.47]{KLMWfr}.

\begin{proposition}\label{prop:quadratic-characterization}
Let $X$ be a Banach space and let $k\geq2$ be an integer. Then $\NA^{(k)}(X,\ell_2^k)\neq\varnothing$ if and only if there exist $f\in S_{X^*}\cap\NA(X,\R)$ and $g=(g_1,\ldots,g_{k-1})\in(X^*)^{k-1}$ such that $f,g_1,\ldots,g_{k-1}$ are linearly independent and $\vartheta(f,g)>0$.
\end{proposition}
\begin{proof}
Assume first that $T\in\NA^{(k)}(X,\ell_2^k)$ with $\|T\|=1$, and choose $x_0\in S_X$ such that $\|Tx_0\|_2=1$. After composing $T$ with a rotation on $\ell_2^k$, we may assume that $Tx_0=(1,0,\ldots,0)$. Write $Tx=(f(x),g_1(x),\ldots,g_{k-1}(x))$ for every $x\in X$, and $g=(g_1,\ldots,g_{k-1})$. Then, $f(x_0)=1$ and $g(x_0)=0$. Since $\|T\|=1$, we have
    $$
    f(x)^2+\|g(x)\|_2^2=\|Tx\|_2^2\leq1,
    $$
for every $x\in B_X$. In particular, $f\in S_{X^*}\cap\NA(X,\R)$. Moreover, $f,g_1,\ldots,g_{k-1}$ are linearly independent since $T$ has rank $k$. Whenever $g(x)\neq0$, we have
$$\frac{1-f(x)^2}{\|g(x)\|_2^2}\geq1.$$
Therefore, $\vartheta(f,g)\geq1$.

Conversely, assume that $f\in S_{X^*}\cap\NA(X,\R)$, that $f,g_1,\ldots,g_{k-1}$ are linearly independent and that $\vartheta(f,g)>0$. Choose $\lambda>0$ such that $\lambda^2\leq\vartheta(f,g)$, and define $T:X\to\ell_2^k$ by
$$Tx:=(f(x),\lambda g_1(x),\ldots,\lambda g_{k-1}(x)).$$
For every $x\in B_X$ with $g(x)\neq0$, we have
$$\lambda^2\leq\frac{1-f(x)^2}{\|g(x)\|_2^2},$$
and hence $f(x)^2+\lambda^2\|g(x)\|_2^2\leq1$. The same inequality is immediate when $g(x)=0$. Thus, $\|Tx\|_2\leq1$ for every $x\in B_X$, and therefore $\|T\|\leq1$.
Since $f$ attains its norm, there exists $x_0\in S_X$ such that $f(x_0)=1$. We claim that $g(x_0)=0$. Otherwise,
$$\vartheta(f,g)\leq\frac{1-f(x_0)^2}{\|g(x_0)\|_2^2}=0,$$
which is impossible. Hence, $Tx_0=(1,0,\ldots,0)$, and therefore $\|Tx_0\|_2=1$. It follows that $\|T\|=1$ and that $T$ attains its norm. Since $f,g_1,\ldots,g_{k-1}$ are linearly independent, we have $\rank T=k$.
\end{proof}

\begin{corollary}\label{cor:quadratic-nonexistence}
Let $X$ be a Banach space and let $k\geq2$ be an integer. The following assertions are equivalent:
\begin{enumerate}
\item[\rm{(a)}] $\NA^{(k)}(X,\ell_2^k)=\varnothing$.
\item[\rm{(b)}] $\vartheta(f,g)=0$ for every $f\in S_{X^*}\cap\NA(X,\R)$ and every $g=(g_1,\ldots,g_{k-1})\in(X^*)^{k-1}$ such that $f,g_1,\ldots,g_{k-1}$ are linearly independent.
\item[\rm{(c)}] For every $f$ and $g$ as in \textup{(b)} and every $\varepsilon>0$, there exists $x\in B_X$ such that $g(x)\neq0$ and $1-f(x)^2<\varepsilon\|g(x)\|_2^2$.
\end{enumerate}
\end{corollary}
\begin{proof}
The equivalence between (a) and (b) follows from Proposition~\ref{prop:quadratic-characterization}, whereas the equivalence between (b) and (c) follows from the definition of $\vartheta(f,g)$.
\end{proof}

Let $(X,p)$ be a Banach space.
For $f\in S_{(X,p)^*}$ and $g=(g_1,\ldots,g_m)\in(X^*)^m$ such that $f,g_1,\ldots,g_m$ are linearly independent, define
$$\Gamma(f,g):=\sup_{t\in\R^m\setminus\{0\}}\frac{p^*\left(f+\sum_{j=1}^mt_jg_j\right)^2-1}{\|t\|_2^2}.$$
For $m=1$, the following proposition is a quantitative form of \cite[Proposition~13.37]{KLMWfr}.

\begin{proposition}\label{prop:theta-dual}
Let $(X,p)$ be a Banach space.
For any $f\in S_{(X,p)^*}$ and $g=(g_1,\ldots,g_m)\in(X^*)^m$ such that $f,g_1,\ldots,g_m$ are linearly independent, we have
$$\vartheta(f,g)=\frac{1}{\Gamma(f,g)},$$
where $1/{+\infty}:=0$.
\end{proposition}

\begin{proof}
For any $\lambda>0$, define $T_\lambda : (X,p)\to \ell_2^{m+1}$ by $T_\lambda x:=(f(x),\lambda g_1(x),\ldots,\lambda g_m(x))\in\ell_2^{m+1}$.
By the definition of $\vartheta(f,g)$, we have that $\|T_\lambda\|\leq1$ if and only if $\lambda^2\leq\vartheta(f,g)$. On the other hand, by duality with the Euclidean norm, $\|T_\lambda\|\leq1$ if and only if
    $$p^*\left(af+\lambda\sum_{j=1}^mb_jg_j\right)\leq\sqrt{a^2+\|b\|_2^2},$$
for all $a\in\R$ and $b\in\R^m$. For $a\neq0$, this means that
    $$
    p^*\left(f+\sum_{j=1}^mt_jg_j\right)^2\leq1+\frac{\|t\|_2^2}{\lambda^2},
    $$
for all $t\in\R^m$. The case $a=0$ follows by replacing $t$ by $st$ and letting $s\to\infty$. Hence, $\|T_\lambda\|\leq1$ if and only if $\Gamma(f,g)\leq\lambda^{-2}$. Taking the supremum over all admissible $\lambda^2$ gives the formula.
\end{proof}

For a seminorm $r$ on $X^*$, $f\in X^*\setminus\{0\}$, $g\in X^*$, and $t\geq 0$, put
    $$
    \Delta_r(f,g;t):=r(f+tg)+r(f-tg)-2r(f),
    $$
and define
    $$
    \mathcal{Q}_r(f):=\left\{g\in X^*\colon \Delta_r(f,g;t)=O(t^2)\text{ as }t\to0^+\right\}.
    $$
By convexity, $\Delta_r(f,g;t)\geq0$.

The condition defining $\mathcal Q_{p^*}(f)$ can be expressed in terms of the classical directional moduli of smoothness. Indeed, given $(X,p)$ a Banach space, recall that the \textit{directional modulus of smoothness of the dual norm $p^*$ at $f\in S_{X^*}$ in the directions of $D\subseteq X^*$} is defined as
    $$
    \rho_{p^*}(t,f,D) := \sup_{g\in D} \left\{ \frac{p^*(f+tg)+p^*(f-tg)}{2}-1 \right\},
    $$
for every $t\geq0$.
Thus, for $D=\{g\}$ with $g\in X^*$, we have
    $$
    \rho_{p^*}(t,f,\{g\}) = \frac{1}{2}\Delta_{p^*}(f,g;t),
    $$
for all $t\geq0$.
Hence, $\mathcal{Q}_{p^*}(f)$ consists of the directions along which the dual norm $p^*$ has smoothness of power type two at $f$.

\begin{proposition}\label{prop:quadratic-directions}
Let $(X,p)$ be a Banach space and $f\in X^*\setminus\{0\}$. Then $\mathcal{Q}_{p^*}(f)$ is a vector subspace of $X^*$ containing $\R f$.
If $f\in S_{(X,p)^*}$ and $k\geq2$ is an integer, then the following assertions are equivalent:
\begin{enumerate}
\item[\rm{(a)}] There exists $g=(g_1,\ldots,g_{k-1})\in(X^*)^{k-1}$ such that $f,g_1,\ldots,g_{k-1}$ are linearly independent and $\vartheta(f,g)>0$.
\item[\rm{(b)}] $\dim\mathcal{Q}_{p^*}(f)\geq k$.
\end{enumerate}
\end{proposition}
\begin{proof}
Scalar multiplication is immediate. For addition, convexity gives
$$\Delta_{p^*}(f,g+h;t)\leq\frac12\Delta_{p^*}(f,g;2t)+\frac12\Delta_{p^*}(f,h;2t),$$
so $\mathcal{Q}_{p^*}(f)$ is a vector subspace. Moreover, for every $\alpha\in\R$ and all sufficiently small $t$,
$$p^*(f+t\alpha f)+p^*(f-t\alpha f)-2p^*(f)=0,$$
and hence $\R f\subseteq\mathcal{Q}_{p^*}(f)$.

\medskip

Assume first (a). By Proposition~\ref{prop:theta-dual}, we have $\Gamma(f,g)<\infty$. Let $u=(u_1,\ldots,u_{k-1})\in\R^{k-1}$ and set $h:=\sum_{j=1}^{k-1}u_jg_j$. 
If $u\neq 0$, by definition of $\Gamma(f,g)$ it yields $p^*(f\pm sh)^2\leq1+\Gamma(f,g)\|u\|_2^2s^2$ for all $s\in\R$. Note that the same inequality holds for $u=0$.
Therefore, for any $s\in \mathbb R$,
    \begin{align*}
        \Delta_{p^*}(f,h;s) & \leq 2\sqrt{1+\Gamma(f,g)\|u\|_2^2s^2}-2 \\
        & = \frac{2\Gamma(f,g)\|u\|_2^2}{\sqrt{1+\Gamma(f,g)\|u\|_2^2s^2}+1}s^2 \\
        & \leq \Gamma(f,g)\|u\|_2^2s^2.
    \end{align*}
Hence, $\Delta_{p^*}(f,h;s)=O(s^2)$.
Thus, $h\in\mathcal{Q}_{p^*}(f)$ which yields $\spn\{f,g_1,\ldots,g_{k-1}\}\subseteq\mathcal{Q}_{p^*}(f)$, and (b) holds.

\medskip

Conversely, assume (b), and choose $h_1,\ldots,h_{k-1}\in\mathcal{Q}_{p^*}(f)$ such that $f,h_1,\ldots,h_{k-1}$ are linearly independent. Fix $1\leq j\leq k-1$ and put
    $$
    \phi_j(t):=p^*(f+th_j).
    $$
The function $\phi_j$ is convex, and the condition $\phi_j(t)+\phi_j(-t)-2\phi_j(0)=O(t^2)$ implies that $\phi_j$ is differentiable at $0$.
Thus, $\phi_j$ lies above its tangent line at $0$. 
Hence both $\phi_j(t)-\phi_j(0)-\phi_j'(0)t$ and $\phi_j(-t)-\phi_j(0)+\phi_j'(0)t$ are non-negative, and their sum is $\phi_j(t)+\phi_j(-t)-2\phi_j(0)=O(t^2)$. Therefore, observe that $\phi_j(t)=\phi_j(0)+\phi_j'(0)t+O(t^2)$, that is, $p^*(f+th_j)=1+\phi_j'(0)t+O(t^2)$. Setting $g_j:=h_j-\phi_j'(0)f$, note that $p^*(f+tg_j)=1+O(t^2)$
as $t\to0^+$. Moreover, $f,g_1,\ldots,g_{k-1}$ are linearly independent. 

Choose $\delta>0$ and $C>0$ such that $p^*(f+\tau g_j)\leq1+C\tau^2$ for all $\tau\in\R$ with $|\tau|\leq\delta$ and all $1\leq j\leq k-1$. If $t=(t_1,\ldots,t_{k-1})\in\R^{k-1}$ satisfies $\|t\|_2\leq \frac{\delta}{k-1}$, then $|(k-1)t_j|\leq\delta$ for every $1\leq j\leq k-1$, and the convexity of $p^*$ gives
    \begin{align*}
        p^*\left(f+\sum_{j=1}^{k-1}t_jg_j\right) & = p^*\left(\frac{1}{k-1}\sum_{j=1}^{k-1}(f+(k-1)t_jg_j)\right) \\
        & \leq \frac{1}{k-1}\sum_{j=1}^{k-1}\left(1+C(k-1)^2t_j^2\right)=1+C(k-1)\|t\|_2^2.
    \end{align*}
Hence the quotient in the definition of $\Gamma(f,g)$ is bounded for $0<\|t\|_2\leq\frac{\delta}{k-1}$. For $\|t\|_2\geq\frac{\delta}{k-1}$, the triangle and Cauchy--Schwarz inequalities give
    $$
    p^*\left(f+\sum_{j=1}^{k-1}t_jg_j\right)\leq1+\left(\sum_{j=1}^{k-1}p^*(g_j)^2\right)^{1/2}\|t\|_2,
    $$
and hence a uniform bound for the same quotient. Thus, $\Gamma(f,g)<\infty$, and Proposition~\ref{prop:theta-dual} yields $\vartheta(f,g)>0$.
\end{proof}

\begin{corollary}\label{cor:dual-quadratic-nonexistence}
Let $(X,p)$ be a Banach space and let $k\geq2$ be an integer. Then, $\NA^{(k)}((X,p),\ell_2^k)\neq\varnothing$ if and only if $\dim\mathcal{Q}_{p^*}(f)\geq k$ for some $f\in S_{(X,p)^*}\cap\NA((X,p),\R)$. In particular, $\NA^{(2)}((X,p),\ell_2^2)=\varnothing$ if and only if $\mathcal{Q}_{p^*}(f)=\R f$ for every $f\in S_{(X,p)^*}\cap\NA((X,p),\R)$.
\end{corollary}
\begin{proof}
Apply Corollary~\ref{cor:quadratic-nonexistence} and Proposition~\ref{prop:quadratic-directions}. Since $\R f\subseteq\mathcal{Q}_{p^*}(f)$, the condition $\dim\mathcal{Q}_{p^*}(f)\geq2$ means that $\mathcal{Q}_{p^*}(f)\neq\R f$.
\end{proof}

Combining Corollary~\ref{cor:dual-quadratic-nonexistence} with the form \eqref{eq:smoothing-dual} of the new dual norm and the preservation of norm-attaining functionals in Lemma~\ref{lem:smoothing}, we obtain the following.

\begin{corollary}\label{cor:dgs-rank-k}
Let $(X,q)$, $G$, $U$ and $p$ be as in Lemma~\ref{lem:smoothing}, and let $k\geq2$ be an integer. Then $\NA^{(k)}((X,p),\ell_2^k)\neq\varnothing$ if and only if $\NA^{(k)}((X,q),\ell_2^k)\neq\varnothing$.
\end{corollary}
\begin{proof}
Put $\rho(f):=\|U^*f\|_G$ for $f\in X^*$. By \eqref{eq:smoothing-dual}, $p^*=q^*+\rho$, and hence
$$\Delta_{p^*}(f,g;t)=\Delta_{q^*}(f,g;t)+\Delta_\rho(f,g;t)$$
for all $f\in X^*\setminus\{0\}$, $g\in X^*$ and $t\in\R$. Both summands are non-negative, so $g\in\mathcal{Q}_{p^*}(f)$ if and only if $g\in\mathcal{Q}_{q^*}(f)$ and $\Delta_\rho(f,g;t)=O(t^2)$.

Since $U$ has dense range, $U^*$ is injective. Let $f\in X^*\setminus\{0\}$ and $g\in X^*$. Then $U^*f\neq0$, and
$$\rho(f+tg)=\left(\|U^*f\|_G^2+2t\ip{U^*f}{U^*g}_G+t^2\|U^*g\|_G^2\right)^{1/2}$$
is a twice continuously differentiable function of $t$ near $0$. By Taylor's theorem, $\Delta_\rho(f,g;t)=O(t^2)$ as $t\to0^+$. Therefore,
$$\mathcal{Q}_{p^*}(f)=\mathcal{Q}_{q^*}(f)\quad\text{for every}\quad f\in X^*\setminus\{0\}.$$

For $c>0$, we have $\Delta_{p^*}(cf,g;t)=c\,\Delta_{p^*}(f,g;t/c)$, so $\mathcal{Q}_{p^*}(cf)=\mathcal{Q}_{p^*}(f)$, and similarly for $q^*$. By Lemma~\ref{lem:smoothing}, $\NA((X,p),\R)=\NA((X,q),\R)$. If $f\neq0$ belongs to the latter sets, then $f/p^*(f)\in S_{(X,p)^*}\cap\NA((X,p),\R)$, $f/q^*(f)\in S_{(X,q)^*}\cap\NA((X,q),\R)$, and
$$\mathcal{Q}_{p^*}\left(\frac{f}{p^*(f)}\right)=\mathcal{Q}_{p^*}(f)=\mathcal{Q}_{q^*}(f)=\mathcal{Q}_{q^*}\left(\frac{f}{q^*(f)}\right).$$
Hence, the condition in Corollary~\ref{cor:dual-quadratic-nonexistence} holds for $(X,p)$ if and only if it holds for $(X,q)$.
\end{proof}

\begin{proposition}\label{prop:one-ray-elimination}
Let $(X,q)$ be a separable Banach space, $f_0\in S_{(X,q)^*}\cap\NA((X,q),\R)$, and $\varepsilon>0$. Then, there exist a compact operator $S:\ell_2\to X$ and an equivalent norm $p$ on $X$ such that
$$B_p=B_q+\varepsilon S(B_{\ell_2}),$$
and the following assertions hold:
\begin{enumerate}
\item[\rm{(a)}] $\overline{S(\ell_2)}=\ker f_0$;
\item[\rm{(b)}] $p^*(f)=q^*(f)+\varepsilon\|S^*f\|_2$ for every $f\in X^*$;
\item[\rm{(c)}] $\NA((X,p),\R)=\NA((X,q),\R)$ as subsets of $X^*$;
\item[\rm{(d)}] for $f\in X^*\setminus\{0\}$,
    $$
    \mathcal{Q}_{p^*}(f)=
    \begin{cases}
    \mathcal{Q}_{q^*}(f) & \text{if }f\notin\R f_0,\\ 
    \R f_0 & \text{if }f\in\R f_0.
    \end{cases}
    $$
\end{enumerate}
\end{proposition}

\begin{proof}
Choose a sequence $(u_n)_n\subseteq B_q\cap\ker f_0$ whose linear span is dense in $\ker f_0$, and define $S:\ell_2\to X$ by
    $$
    S(a):=\sum_{n=1}^\infty2^{-n}a_nu_n.
    $$
Note that $S$ is compact and, since every $u_n\in S(\ell_2)$, we have $\overline{S(\ell_2)}=\ker f_0$, which proves (a). Define $p$ by
    $$
    B_p:=B_q+\varepsilon S(B_{\ell_2}).
    $$
As $S(B_{\ell_2})$ is compact, $B_p$ is the closed unit ball of an equivalent norm. The support function of the Minkowski sum gives $p^*(f)=q^*(f)+\varepsilon\|S^*f\|_2$, proving (b).

Since every functional attains its supremum on the compact set $S(B_{\ell_2})$, a functional which attains its $q^*$-norm also attains its $p^*$-norm. Conversely, assume that $f$ attains its $p^*$-norm at a point $u+\varepsilon Sa$, where $u\in B_q$ and $a\in B_{\ell_2}$. Then,
    $$
    f(u)+\varepsilon(S^*f)(a)=q^*(f)+\varepsilon\|S^*f\|_2.
    $$
Both terms on the left are bounded above by the corresponding terms on the right, so equality holds in both. In particular, $f(u)=q^*(f)$. This proves (c).

Finally, put $\rho(f):=\|S^*f\|_2$. By (a) we have $\ker S^*=S(\ell_2)^\perp=\R f_0$, so $S^*f\neq0$ whenever $f\notin\R f_0$. At such $f$ the Hilbert norm is twice differentiable, and therefore $\Delta_\rho(f,g;t)=O(t^2)$ for every fixed $g\in X^*$. Since $\Delta_{p^*}(f,g;t)=\Delta_{q^*}(f,g;t)+\varepsilon\Delta_\rho(f,g;t)$ and both summands are non-negative, it follows that $\mathcal{Q}_{p^*}(f)=\mathcal{Q}_{q^*}(f)$ for every $f\notin\R f_0$. If $f\in\R f_0\setminus\{0\}$, then $S^*f=0$ and $\Delta_\rho(f,g;t)=2|t|\|S^*g\|_2$. Thus $g\in\mathcal{Q}_{p^*}(f)$ only if $S^*g=0$, equivalently $g\in\R f_0$. The reverse inclusion is automatic. This proves (d).
\end{proof}

\begin{remark}\label{rem:one-ray}
Proposition~\ref{prop:one-ray-elimination} shows that a compact Minkowski enlargement can eliminate all rank-two quadratic directions attached to one prescribed norm-attaining ray without changing the set of norm-attaining functionals.
\end{remark}

\section{Mart\'in's space}\label{sec:martin}
We use the notation $\NA^{(2)}(X,\ell_2^2)$ and $\vartheta(f,g)$ of Section~\ref{sec:rank-k}.

\subsection{Norm-attaining operators of finite rank}
We now give an affirmative answer to \cite[Open Problem~7]{Mar}. In fact, Theorem~\ref{thm:finite-short} yields norm-attaining operators of every prescribed finite rank allowed by the dimension of the range space.

\begin{theorem}\label{thm:finite-short}
Let $\X$ be Mart\'in's space as in Subsection~\ref{subsec:space}. For every Banach space $F$ and every positive integer $k\le\dim F$, the set $\NA(\X,F)$ contains an operator of rank $k$.
\end{theorem}
\begin{proof}
Let $(e_n)_n$ be the canonical basis for $\ell_1$. Choose linearly independent $u_1,\ldots,u_k\in F$ and define $S\in\Lin((c_0,\tn{\cdot}),F)$ by
    $$
    S=\sum_{n=1}^k e_n\otimes u_n.
    $$
By (M1) in Subsection~\ref{subsec:space}, $\NA((c_0,\tn{\cdot}),\R)$ is the set of finitely supported sequences.
Hence Lemma~\ref{lem:densesubspace} shows that $S$ attains its norm.

Set $q(z)=\tn{z}$ and $s(z)=\|Rz\|_W$. By (M5), $s$ is a continuous seminorm with $s\le2q$. Since the norm $p$ of $\X$ satisfies $p=q+s$, Lemma~\ref{lem:lift} gives $A\in\NA(\X,F)$ with $\ran A=\ran S$. Consequently, $\rank A=k$.
\end{proof}

The density assertion of Lemma~\ref{lem:densesubspace} cannot be applied directly to $\X$: by \cite{Mar}, $\NA(\X,\R)$ contains no cone generated by two linearly independent functionals. In particular, it contains no dense linear subspace of $\X^*$.

We note that Lemma~\ref{lem:lift}, which underlies Theorem~\ref{thm:finite-short}, holds for an arbitrary continuous seminorm, and that its argument also applies to absolute combinations of a norm and a seminorm. Recall that a norm $\nu$ on $\R^2$ is \emph{absolute} if $\nu(\xi,\eta)=\nu(|\xi|,|\eta|)$ for all $\xi,\eta\in\R$, and \emph{normalized} if $\nu(1,0)=\nu(0,1)=1$; we write $\nu^*$ for its dual norm. An absolute norm is non-decreasing in each variable on $[0,\infty)^2$, so $p_\nu$ below is a norm.

\begin{proposition}\label{prop:absolute-lifting}
Let $q$ be a norm on a vector space $X$, let $s$ be a continuous seminorm on $(X,q)$, and let $\nu$ be an absolute normalized norm on $\R^2$. Define
$$p_\nu(x):=\nu(q(x),s(x))$$
for every $x\in X$. Assume that there are $T_0\in\NA^{(2)}((X,q),\ell_2^2)$ and $x_0\in X$ such that $q(x_0)=\|T_0\|=\|T_0x_0\|_2=1$, and assume that there exists a norming functional $(\alpha,\beta)$ for $(1,s(x_0))$ with respect to $\nu$ satisfying $\alpha>0$, $\beta\geq0$, $\nu^*(\alpha,\beta)=1$, and $\alpha+\beta s(x_0)=\nu(1,s(x_0))$. Then, $\NA^{(2)}((X,p_\nu),\ell_2^2)\neq\varnothing$.
\end{proposition}
\begin{proof}
After composing $T_0$ with a rotation on $\ell_2^2$, write $T_0x=(f_0(x),g_0(x))$ so that $T_0x_0=(1,0)$. Define a linear functional $h_0$ on $\R x_0$ by $h_0(\lambda x_0):=\lambda s(x_0)$ for every $\lambda\in\R$. Since $s$ is a seminorm, we have $|h_0(\lambda x_0)|=|\lambda|s(x_0)=s(\lambda x_0)$.
Therefore, by the Hahn--Banach theorem for seminorms, there exists $h\in X^*$ such that $h(x_0)=s(x_0)$ and $|h(x)|\leq s(x)$ for every $x\in X$. Define
$$Tx:=\left(\alpha f_0(x)+\beta h(x),\alpha g_0(x)\right).$$
For every $x\in X$, we have
$$\|Tx\|_2\leq\alpha\|(f_0(x),g_0(x))\|_2+\beta|h(x)|\leq\alpha q(x)+\beta s(x).$$
Since $\nu^*(\alpha,\beta)=1$, duality for the norm $\nu$ gives
$$\alpha q(x)+\beta s(x)\leq\nu(q(x),s(x))=p_\nu(x).$$
Hence, $\|T\|\leq1$.

At $x_0$, we have $Tx_0=(\alpha+\beta s(x_0),0)$, and therefore
$$\|Tx_0\|_2=\alpha+\beta s(x_0)=\nu(1,s(x_0))=p_\nu(x_0).$$
It follows that
$$\left\|T\left(\frac{x_0}{p_\nu(x_0)}\right)\right\|_2=1.$$
Thus, $\|T\|=1$ and $T$ attains its norm.

It remains to prove that $T$ has rank two. Assume, by means of contradiction, that $T$ has rank one. Then, the two coordinate functionals $\alpha f_0+\beta h$ and $\alpha g_0$ are linearly dependent. Hence, there exists $\lambda\in\R$ such that $\alpha g_0=\lambda(\alpha f_0+\beta h)$. Evaluating at $x_0$, we obtain $0=\lambda(\alpha+\beta s(x_0))=\lambda\nu(1,s(x_0))$. Since $\nu(1,s(x_0))>0$, it follows that $\lambda=0$. As $\alpha>0$, we then obtain $g_0=0$, contradicting the fact that $T_0$ has rank two. Therefore, $T\in\NA^{(2)}((X,p_\nu),\ell_2^2)$, which finishes the proof.
\end{proof}

\begin{corollary}\label{cor:ell-r-lifting}
Let $q$ and $s$ be as in Proposition~\ref{prop:absolute-lifting}, let $1\leq r<\infty$ and define
$$p_r(x):=\left(q(x)^r+s(x)^r\right)^{1/r}.$$
If $\NA^{(2)}((X,q),\ell_2^2)\neq\varnothing$, then $\NA^{(2)}((X,p_r),\ell_2^2)\neq\varnothing$.
\end{corollary}
\begin{proof}
Apply Proposition~\ref{prop:absolute-lifting} to the absolute norm $\nu_r(\xi,\eta):=\left(|\xi|^r+|\eta|^r\right)^{1/r}$. At every point $(1,b)$ with $b\geq0$, this norm admits a norming functional $(\alpha,\beta)$ with $\alpha>0$.
\end{proof}

\subsection{Mart\'in's base space}\label{subsec:base-space}
We next show that the obstruction is already abundant in every Banach space whose norm-attaining functionals contain a norm-dense linear subspace. For a Banach space $X$, define
$$\mathcal{H}(X):=\left\{f\in S_{X^*}\cap\NA(X,\R) \colon \text{there exists }g\in X^*\setminus\R f\text{ such that }\vartheta(f,g)>0\right\}.$$
Thus, $\mathcal{H}(X)$ consists of the functionals which can occur as the first coordinate of a rotated rank-two norm-attaining operator into $\ell_2^2$.

\begin{proposition}\label{prop:density-bad-coordinates}
Let $X$ be a Banach space. Assume that $\NA(X,\R)$ contains a norm-dense linear subspace $D\subseteq X^*$ with $\dim D\geq2$. Then, $\mathcal{H}(X)$ is norm-dense in $S_{X^*}$.
\end{proposition}
\begin{proof}
Note that $D\cap S_{X^*}$ is norm-dense in $S_{X^*}$. Thus, it is enough to approximate an arbitrary element $f\in D\cap S_{X^*}$ by elements of $\mathcal{H}(X)$.

Choose $h\in D$ linearly independent of $f$. For every $s>0$, define $T_s:X\to\ell_2^2$ by
$$T_sx:=(f(x),sh(x)).$$
Since $f$ and $h$ are linearly independent,
$$[\ker T_s]^\perp=\spn\{f,h\}\subseteq D\subseteq\NA(X,\R).$$
By Lemma~\ref{lem:densesubspace}(a), $T_s\in\NA(X,\ell_2^2)$. Moreover, $T_s$ has rank two.

Since $\|f\|=1$, we have $\|T_s\|\geq1$. On the other hand, $\|T_sx\|_2\leq\sqrt{1+s^2\|h\|^2}\,\|x\|$ for every $x\in X$. Hence,
$$1\leq\|T_s\|\leq\sqrt{1+s^2\|h\|^2},$$
and therefore $\|T_s\|\longrightarrow1$ as $s\to0^+$.

Choose $x_s\in S_X$ such that $\|T_sx_s\|_2=\|T_s\|$. Replacing $x_s$ by $-x_s$, if necessary, we may assume that $f(x_s)\geq0$. Denote $a_s:=\frac{f(x_s)}{\|T_s\|}$ and $b_s:=\frac{sh(x_s)}{\|T_s\|}$. Then, $a_s^2+b_s^2=1$. Moreover,
    $$
    |b_s|\leq\frac{s\|h\|}{\|T_s\|}\longrightarrow0.
    $$
Since $a_s\geq0$, it follows that $a_s\longrightarrow1$.

Consider the orthogonal operator $\pi_s:\ell_2^2\to\ell_2^2$ given by
$$\pi_s(u,v):=(a_su+b_sv,-b_su+a_sv).$$
Then, $\pi_s(a_s,b_s)=(1,0)$. Define the normalized rotated operator $\widetilde T_s:=\pi_s\circ\frac{T_s}{\|T_s\|}$. Writing $\widetilde T_sx=(f_s(x),g_s(x))$, we have
$$f_s=\frac{a_sf+b_ssh}{\|T_s\|}\quad\text{and}\quad g_s=\frac{-b_sf+a_ssh}{\|T_s\|}.$$
The operator $\widetilde T_s$ has norm one, rank two, and satisfies $\widetilde T_sx_s=(1,0)$. Thus, $f_s\in S_{X^*}\cap\NA(X,\R)$, $g_s\notin\R f_s$ and $f_s(x)^2+g_s(x)^2\leq1$ for every $x\in B_X$. It follows that $\vartheta(f_s,g_s)\geq1$, so $f_s\in\mathcal{H}(X)$.

Finally,
$$\|f_s-f\|\leq\left|\frac{a_s}{\|T_s\|}-1\right|\|f\|+\frac{|b_s|s}{\|T_s\|}\|h\|.$$
Since $a_s\longrightarrow1$, $\|T_s\|\longrightarrow1$ and $b_s\longrightarrow0$, we obtain $\|f_s-f\|\longrightarrow0$. Thus, every element of $D\cap S_{X^*}$ belongs to $\overline{\mathcal{H}(X)}$. Since $D\cap S_{X^*}$ is norm-dense in $S_{X^*}$, the result follows.
\end{proof}

\begin{corollary}\label{cor:dense-bad-coordinates}
Let $E=(c_0,\tn{\cdot})$ be the base space in \textup{(M1)}. Then, $\mathcal{H}(E)$ is norm-dense in $S_{E^*}$.
\end{corollary}
\begin{proof}
By (M1), $\NA(E,\R)=\NA(c_0,\R)=c_{00}$. The space $c_{00}$ is a norm-dense linear subspace of $E^*$ and has dimension greater than one. The conclusion follows from Proposition~\ref{prop:density-bad-coordinates}.
\end{proof}

The preceding density has a direct consequence for any attempt to repair all such functionals by enlarging the primal unit ball.

\begin{corollary}\label{cor:no-simultaneous-repair-martin}
Let $E=(c_0,\tn{\cdot})$ be the base space in \textup{(M1)}. There is no closed absolutely convex set $C\subseteq E$ with $B_E\subsetneq C$ such that $\sup_{x\in C}|f(x)|=1$ for every $f\in\mathcal{H}(E)$.
\end{corollary}

\begin{proof}
By Corollary~\ref{cor:dense-bad-coordinates}, $\mathcal H(E)$ is norm-dense in $S_{E^*}$. Assume, by means of contradiction, that such a set $C$ exists. Choose $x_0\in C\setminus B_E$, so that $\|x_0\|>1$. By the Hahn--Banach theorem, there exists $f\in S_{E^*}$ such that $f(x_0)=\|x_0\|>1$. Since $\mathcal H(E)$ is norm-dense in $S_{E^*}$, we may choose $g\in\mathcal H(E)$ sufficiently close to $f$ so that $g(x_0)>1$.
However, this contradicts $\sup_{x\in C}|g(x)|=1$.
\end{proof}

\subsection{Non-norm-attaining operators and the density problem}\label{subsec:nonattaining}
Proposition~\ref{prop:lineable} below provides, for every nonzero Banach space $F$, a norm-dense linear subspace of $\K(\X,F)$ whose nonzero elements fail to attain their norms.
This does not exclude density of $\NA(\X,F)$, and whether $\X$ yields a counterexample to Lindenstrauss property B for a finite-dimensional range remains open.

\begin{proposition}\label{prop:lineable}
For every nontrivial Banach space $F$, the linear subspace 
$$\mathcal M_F:=Y\otimes F=\left\{\sum_{j=1}^k h_j\otimes u_j~:~k<\infty,\ h_j\in Y,\ u_j\in F\right\}$$ 
satisfies $\mathcal M_F\cap \NA(\X,F)=\{0\}$ and $\overline{\mathcal M_F}=\K(\X,F)$.
\end{proposition}
\begin{proof}
We first show that $\NA(\X,\R)\cap Y=\{0\}$ using (M1), (M6) and (M7) of Subsection~\ref{subsec:space}. Assume, by means of contradiction, that $f\in\NA(\X,\R)\cap Y\setminus\{0\}$. Since $Y$ is a linear subspace, we assume that $\|f\|_{\X^*}=1$, and take $x\in S_\X$ with $f(x)=1$.
By (M7), there exist $f_0\in S_{(c_0,\tn{\cdot})^*}$ and $w\in W^*$ such that $f=f_0+R^*w$ and $f_0(x)=\tn{x}$. Thus $f_0$ attains its norm on $(c_0,\tn{\cdot})$ at $\tfrac{x}{\tn{x}}$. On the other hand, (M6) gives $R^*w\in Y$, hence $f_0=f-R^*w\in Y$. This contradicts (M1), since $f_0\ne0$.

Now assume that $A\in\mathcal M_F$ attains its norm at $x\in S_\X$. Then there exists $u^*\in S_{F^*}$ such that $u^*(Ax)=\|A\|$. Therefore $A^*u^*\in Y\cap\NA(\X,\R)$ and thus $A=0$.

For the density assertion, since $Y$ is norm dense in $\X^*$, approximating the finitely many functionals in a finite-rank representation shows that the closure of $\mathcal M_F$ contains every approximable operator.
Moreover, every compact operator from a renorming of $c_0$ is approximable by Lemma~\ref{lem:projections}(a), so we have $\overline{\mathcal M_F}=\K(\X,F)$.
\end{proof}

\subsection*{Acknowledgements}
Manwook Han was supported by Global - Learning \& Academic research institution for Master’s·PhD students, and Postdocs (LAMP) Program of the National Research Foundation of Korea (NRF) grant funded by the Ministry of Education (RS-2024-00445180). He also thanks the Polytechnic University of Madrid for its hospitality during his visit in November 2026. 
Miguel Martín was supported by MICIU/AEI/10.13039/501100011033 and ERDF/EU through the grants PID2021-122126NB-C31 and PID2025-167660NB-I00, and by ``Maria de Maeztu'' Excellence Unit IMAG, funded by MICIU/AEI/10.13039/501100011033 with reference CEX2020-001105-M.
Daniel L. Rodríguez-Vidanes was supported by the Spanish Ministry of Science, Innovation and Universities and the Spanish State Research Agency, project PID2025-171414NA-I00 (MICIU/AEI/10.13039/501100011033).

\subsection*{AI Declaration}
The authors used GPT-6 Astra (OpenAI) to assist with language editing and manuscript organization.
The rank-two case of Theorem~\ref{thm:finite-short}, which answers the question raised in \cite{Mar}, was obtained independently by two of the authors. In June 2026, D.~L.~Rodr\'iguez-Vidanes found a solution without the use of any AI tool. In August 2026, M.~Han obtained a proof of this case with the assistance of Claude Fable 5 (Anthropic). A solution to the same question generated by GPT-5.5 (OpenAI) was also posted in July 2026 on the website \emph{Mathematical Discoveries in the Wild} (\url{https://mathematical-discoveries-in-the-wild.github.io/}). The extension to arbitrary finite rank in Theorem~\ref{thm:finite-short} was obtained by the authors without the use of AI tools.
The authors take full responsibility for the content of this paper.

\end{document}